\documentclass[cmp]{svjour}
\usepackage{amsmath}%
\usepackage{amsfonts}%
\usepackage{amssymb}%
\usepackage{graphicx}

\begin{document}

\title{\textbf{Quantum enveloping algebras with von Neumann regular
Cartan-like generators and the Pierce decomposition}}

\author{Steven Duplij\inst{1}\and Sergey Sinel'shchikov\inst{2}}

\institute{Institut f\"ur Theoretische Physik, Universitat zu
K\"oln, Z\"ulpicher Str. 77, 50937 K\"oln, Germany\thanks{\emph{On
leave of absence from}: Theory Group, Nuclear Physics Laboratory, V.
N. Karazin Kharkov National University, Svoboda Sq. 4, Kharkov
61077, Ukraine,
\texttt{sduplij@gmail.com},
http://webusers.physics.umn.edu/\~{}duplij.},
\email{duplij@thp.uni-koeln.de}\and Mathematics Division, B. I.
Verkin Institute for Low Temperature Physics and Engineering, 47 Lenin
Ave, National Academy of Sciences of Ukraine, Kharkov 61103,
Ukraine, \email{sinelshchikov@ilt.kharkov.ua}}

\date{March 23, 2008}
\dedication{ Dedicated to the memory of our colleague Leonid L.
Vaksman (1951--2007)}
\titlerunning{Quantum enveloping algebras and the Pierce decomposition}

\maketitle

\begin{abstract}
Quantum bialgebras derivable from $U_q(sl_2)$ which contain
idempotents and von Neumann regular Cartan-like generators are
introduced and investigated. Various types of antipodes (invertible
and von Neumann regular) on these bialgebras are constructed, which
leads to a Hopf algebra structure and a von Neumann-Hopf algebra
structure, respectively. For them, explicit forms of some particular
$R$-matrices (also, invertible and von Neumann regular) are
presented, and the latter respects the Pierce decomposition.
\end{abstract}

\section{Introduction}

The language of Hopf algebras \cite{abe,sweedler} is among the
principal tools of studying subjects associated to noncommutative
spaces \cite{connes,madore} and superspaces
\cite{boe/gra/nie,gra/var/fig,sei/wit} appearing as quantization of
commutative ones \cite{wes/bag,gat/gri/roc/sie}. An important
feature of supersymmetric algebraic structures is that their
underlying algebras normally contain idempotents and other zero
divisors \cite{berezin,dup18,rab2}. Therefore, it is reasonable to
render idempotents to some quantum algebras, to study their
properties and the associated Pierce decompositions \cite{pierce}.

In this paper we introduce a new quantum algebra which admits an
embedding of $U_{q}\left(  sl_{2}\right)  $ \cite{dri2,jantzen}.
After adding some extra relations we obtain two worthwhile
algebras that contain idempotents and von Neumann regular
Cartan-like generators. One of the algebras has the Pierce
decomposition which reduces to a direct sum of two ideals and can be
treated as an extended version of the algebra with von Neumann
regular antipode considered in \cite{dup/li3,dup/li2}, while another
one appears to be a Hopf algebra in the sense of the standard
definition \cite{abe}. We distinguish some special cases for which
$R$-matrices of simple form are available. This way both invertible
and von Neumann regular $R$-matrices have been produced, the latter
respecting the Pierce decomposition.

\section{Preliminaries}

We start with recalling briefly some necessary notations and
principal facts about Hopf algebras \cite{abe,cha/pre}. In our
context an algebra $U^{\left( alg\right)  }$ over $\mathbb{C}$ is a
4-tuple $\left(  \mathbb{C},A,\mu ,\eta\right)  $, where $A$ is a
vector space, $\mu:A\otimes A\rightarrow A$ is a multiplication
(alternatively denoted as $\mu\left(  a\otimes b\right)
=a\cdot b$), $\eta:\mathbb{C}\rightarrow A$ is a unit so that $\mathbf{1}%
\overset{def}{=}\eta\left(  1\right)  $, $\mathbf{1\in}A$,
$1\in\mathbb{C}$. The multiplication is assumed to be associative
$\mu\circ\left(  \mu \otimes\mathsf{id}\right)  =\mu\circ\left(
\mathsf{id}\otimes\mu\right)  $ and the unit is characterized by the
property $\mu\circ\left(  \eta \otimes\mathsf{id}\right)
=\mu\circ\left(  \mathsf{id}\otimes\eta\right) =\mathsf{id}$. An
algebra map is a linear map $\psi:U_{1}^{\left(  alg\right)
}\rightarrow U_{2}^{\left(  alg\right)  }$ subject to $\psi\circ\mu_{1}%
=\mu_{2}\circ\left(  \psi\otimes\psi\right)  $ and
$\psi\circ\eta_{1}=\eta _{2}$. A coalgebra $U^{\left(  coalg\right)
}$ is a 4-tuple $\left( \mathbb{C},C,\Delta,\epsilon\right)  $,
where $C$ is an underlying vector space, $\Delta:C\rightarrow
C\otimes C$ is a comultiplication with $\Delta\left(  A\right)
=\sum_{i}\left(  A_{\left(  1\right)  }^{i}\otimes A_{\left(
2\right)  }^{i}\right)  $ in the Sweedler notation, $\epsilon
:C\rightarrow\mathbb{C}$ is a counit. These linear maps are subject
to the
following properties: coassociativity $\left(  \Delta\otimes\mathsf{id}%
\right)  \circ\Delta=\left(  \mathsf{id}\otimes\Delta\right)
\circ\Delta$, the counit property $\left(
\epsilon\otimes\mathsf{id}\right)  \circ \Delta=\left(
\mathsf{id}\otimes\epsilon\right)  \circ\Delta=\mathsf{id}$. A
coalgebra map is a linear map $\varphi:U_{1}^{\left(  coalg\right)
}\rightarrow U_{2}^{\left(  coalg\right)  }$ such that $\left(
\varphi \otimes\varphi\right)
\circ\Delta_{1}=\Delta_{2}\circ\varphi$ and
$\epsilon_{1}=\epsilon_{2}\circ\varphi$. A bialgebra $U^{\left(
bialg\right) }$ is a 6-tuple $\left(
\mathbb{C},B,\mu,\eta,\Delta,\epsilon\right)$ which is an algebra and
coalgebra simultaneously, with the compatibility conditions as follows: $\Delta\circ
\mu=\left(  \mu\otimes\mu\right)  \circ(\mathsf{id\otimes\tau\otimes}\mathsf{id})\circ(\Delta\otimes\Delta$), $\Delta\left(\mathbf{1}\right)=\mathbf{1}\otimes\mathbf{1}$, $\epsilon\circ\mu=\mu_\mathbb{C}\circ(\epsilon\otimes\epsilon)$,
 $\epsilon\left(  \mathbf{1}\right)  =1$; here $\tau$ is the flip of tensor
multiples, $\mu_{\mathbb{C}}$ is the multiplication in the
ground field. A Hopf
algebra $U^{\left( Hopf\right)  }$ is a bialgebra equipped with
antipode, an antimorphism of algebra subject to the relation $\left(
S\otimes\mathsf{id}\right) \circ\Delta=\left(  \mathsf{id}\otimes
S\right)  \circ\Delta=\eta\circ \epsilon$.

Let $q\in\mathbb{C}$ and $q\neq\pm1$,$0$. We start with a definition
of quantum universal enveloping algebra $U_{q}\left(  sl_{2}\right)
$ \cite{dri0}. This is a unital associative algebra $U_{q}^{\left(
alg\right) }\left(  sl_{2}\right)  $ determined by its (Chevalley)
generators $k$,
$k^{-1}$, $e$, $f$, and the relations%
\begin{align}
& k^{-1}k    =\mathbf{1},\ \ \ \ kk^{-1}=\mathbf{1},\label{kk1}\\
& ke    =q^{2}ek,\ \ \ kf=q^{-2}fk,\\
& ef-fe   =\dfrac{k-k^{-1}}{q-q^{-1}}.
\end{align}
The standard Hopf algebra structure on $U_{q}^{(Hopf)}\left(
sl_{2}\right)  $
is determined by%
\begin{align}
\Delta_{0}\left(  k\right)   &  =k\otimes k\label{k0}\\
\Delta_{0}\left(  e\right)   &  =\mathbf{1}\otimes e+e\otimes
k,\;\Delta_{0}\left( f\right)
=f\otimes\mathbf{1}+k^{-1}\otimes f,\label{def}\\
\mathsf{S_{0}}\left(  k\right)   &  =k^{-1},\;\mathsf{S_{0}}\left(  e\right)
=-ek^{-1},\;\mathsf{S_{0}}\left(  f\right)  =-kf,\label{ss}\\
\varepsilon_0\left(  k\right)   &  =1,\ \ \ \varepsilon_0\left(  e\right)
=\varepsilon_0\left(  f\right)  =0. \label{ef}%
\end{align}
The algebra $U_{q}^{(alg)}\left(  sl_{2}\right)  $ is a domain, i.e.
it has no zero divisors and, in particular, no idempotents
\cite{con/kac,jos/let}. A basis of the vector space $U_{q}\left(
sl_{2}\right)  $ is given by the monomials $k^{s}e^{m}f^{n}$, where
$m,n\geq0$ \cite{jantzen}. We denote the Cartan subalgebra of
$U_{q}\left(  sl_{2}\right)  $ by $\mathcal{H}_{0}\left(
\mathbf{1},k,k^{-1}\right)  $.

Our goal is to apply the Pierce decomposition to a suitably extended
version of $U_{q}\left(  sl_{2}\right)  $. It is well known that
there exists one-to-one correspondence between the central
decompositions of unity on idempotents and decompositions of a
module into a direct sum. Therefore we start with generalizing the
Cartan subalgebra in $U_{q}\left(  sl_{2}\right) $ towards von
Neumann regularity property \cite{nashed,rao/mit,cam/mey}.

\section{From the standard \boldmath $U_{q}\left(  sl_{2}\right)  $ to $U_{K+L}$}

Let us consider the generators $K$, $\overline{K}$ satisfying the relations%
\begin{equation}
K\overline{K}K=K,\ \ \ \ \overline{K}K\overline{K}=\overline{K}, \label{k}%
\end{equation}
which are normally referred to as von Neumann regularity
\cite{nashed}. Under
the assumption of commutativity%
\begin{equation}
K\overline{K}=\overline{K}K \label{kk}%
\end{equation}
we have an idempotent $P\overset{def}{=}K\overline{K}=\overline{K}K$
subject
to%
\begin{align}
PK  &  =KP=K,\label{kp}\\
P^{2}  &  =P. \label{p2}%
\end{align}
The commutative algebra generated by $K$, $\overline{K}$ is not
unital (we denote it by $\mathcal{H}\left(  K,\overline{K}\right)
$), because unlike $U_{q}\left(  sl_{2}\right)  $ its relations do
not anticipate unit explicitly, as in (\ref{kk1}). Note that
$\mathcal{H}\left(  K,\overline {K}\right)  $ was considered as a
Cartan-like part of the analog of quantum
enveloping algebra with von Neumann regular antipode $U_{q}^{v}=\mathfrak{vsl}%
_{q}\left(  2\right)  $ introduced by Duplij and Li
\cite{dup/li3,dup/li2}. The associated unital algebra derived by an
exterior attachment of unit
$\mathcal{H}\left(  \mathbf{1},K,\overline{K}\right)  \overset{def}{=}\mathcal{H}%
\left(  K,\overline{K}\right)  \oplus\mathbb{C}\mathbf{1}$ also
appears in \cite{dup/li3,dup/li2} as a part of
$U_{q}^{w}=\mathfrak{wsl}_{q}\left( 2\right)  $.

Observe that $\mathcal{H}\left(  \mathbf{1},K,\overline{K}\right)  $
contains one more idempotent $\left(  \mathbf{1}-P\right)
^{2}=\left( \mathbf{1}-P\right)  $. Therefore, we
introduce another copy of the same algebra (we denote it by $\mathcal{H}%
\left(  L,\overline{L}\right)  $) with generators $L$ and
$\overline{L}$
subject to similar relations as for $K$, $\overline{K}$ above%
\begin{equation}
L\overline{L}L-L=0,\ \ \ \ \overline{L}L\overline{L}-\overline{L}=0. \label{l}%
\end{equation}
Under the commutativity assumption
\begin{equation}
L\overline{L}=\overline{L}L \label{ll}%
\end{equation}
the idempotent $Q\overset{def}{=}L\overline{L}=\overline{L}L$ satisfies%
\begin{align}
QL  &  =LQ=L,\label{lq}\\
Q^{2}  &  =Q. \label{q2}%
\end{align}
If there are no additional relations between $K,\overline{K}$ and
$L,\overline{L}$, the nonunital algebras $\mathcal{H}\left(
K,\overline {K}\right)  $ and $\mathcal{H}\left(
L,\overline{L}\right)  $ can form a free product only. On the other
hand we merge together the unital algebras $\mathcal{H}\left(
\mathbf{1},K,\overline{K}\right)  $ and $\mathcal{H}\left(
\mathbf{1},L,\overline{L}\right)  $ so that their units are
identified and add one more relation, the decomposition of
unity%
\begin{equation}
P+Q=\mathbf{1} \label{pq1}%
\end{equation}
in order to produce the Pierce decomposition \cite{pierce} of the
resulting algebra $\mathcal{H}\left(
\mathbf{1},K,\overline{K},L,\overline{L}\right)  $, which reduces to
the direct product since $QP=PQ=0$.

It follows from (\ref{kp}), (\ref{lq}) and (\ref{pq1}) that%
\begin{equation}
KL=\overline{L}K=LK=K\overline{L}=\overline{K}L=L\overline{K}=0. \label{lk}
\end{equation}
The new (as compared to \cite{dup/li3,dup/li2}) noninvertible
generators $L$, $\overline{L}$ are introduced to justify the
following

\begin{lemma}\label{invelt}
The sum $aK+bL$ is invertible, and its inverse is
$a^{-1}\bar{K}+b^{-1}\bar {L}$ , where $a,b\in\mathbb{R}\diagdown0$.
\end{lemma}

\begin{proof}
reduces to a computation which involves (\ref{pq1}) and (\ref{lk})
as
\begin{equation}
\left(  aK+bL\right)  \left(  a^{-1}\bar{K}+b^{-1}\bar{L}\right)
=K\bar
{K}+L\bar{L}=P+Q=\mathbf{1}. \label{kkll}%
\end{equation}
\end{proof}
This allows us to consider a two-parameter family of morphisms for
the Cartan subalgebra $\mathbf{\Phi}_{\mathcal{H}}^{\left(
a,b\right) }:\mathcal{H}_0\left(
\mathbf{1},k,k^{-1}\right)  \rightarrow\mathcal{H}\left(  \mathbf{1},K,\overline{K}%
,L,\overline{L}\right)  $ given by%
\begin{equation}
k\rightarrow aK+bL,\ \ \ \ \ \ k^{-1}\rightarrow a^{-1}\overline{K}%
+b^{-1}\overline{L}. \label{kab}%
\end{equation}

\begin{proposition}\label{kerf1}
The map $\mathbf{\Phi}_{\mathcal{H}}^{\left( a,b\right)
}$ is an embedding, i.e. $\ker\mathbf{\Phi}_{\mathcal{H}}^{\left(
a,b\right) }=0$.
\end{proposition}

\begin{proof}
Use (\ref{kab}) to define a homomorphism
$\bar{\mathbf{\Phi}}_{\mathcal{H}}^{\left(
a,b\right)  }$ from the free algebra $\mathcal{\bar{H}}_0\left(  \mathbf{1},k,k^{-1}%
\right)  $ generated by $\mathbf{1}$, $k$, $k^{-1}$ into the free
algebra $\mathcal{\bar{H}}\left(
\mathbf{1},K,\overline{K},L,\overline{L}\right)  $ generated by
$\mathbf{1}$, $K$, $\overline{K}$, $L$, $\overline{L}$. We claim
that $\bar{\mathbf{\Phi} }_{\mathcal{H}}^{\left(  a,b\right)  }$ is
an embedding. In fact, if not, then
$\bar{\mathbf{\Phi}}_{\mathcal{H}}^{\left( a,b\right)  }$
annihilates some nonzero element of $\mathcal{\bar{H}}_0\left(
\mathbf{1},k,k^{-1}\right)  $. This element can be treated as a
``noncommutative polynomial'' in three indeterminates $\mathbf{1}$,
$k$, $k^{-1}$. Because the linear change of variables (\ref{kab}) is
non-degenerate, we obtain a nontrivial polynomial in $\mathbf{1}$,
$K$, $\overline{K}$, $L$, $\overline{L}$, which cannot be zero in
the free algebra $\mathcal{\bar{H}}\left(
\mathbf{1},K,\overline{K},L,\overline{L}\right) $. What remains is
to observe that $\mathbf{\Phi}_{\mathcal{H}}^{\left(
a,b\right) }$establishes one-to-one correspondence between the
relations in $\mathcal{H}_0\left(\mathbf{1},k,k^{-1}\right)  $
and those induced on the image of
$\mathbf{\Phi}_{\mathcal{H}}^{\left( a,b\right) }$, which
already implies our statement for the morphism
$\mathbf{\Phi}_{\mathcal{H}}^{\left(  a,b\right) }$ between
the quotient algebras $\mathcal{H}_{0}\left(
\mathbf{1},k,k^{-1}\right) $ and $\mathcal{H}\left(
\mathbf{1},K,\overline{K},L,\overline{L}\right)  $.
\end{proof}

Now we are in a position to add two more generators $E$ and $F$,
along with
additional relations%
\begin{align}
& \left(  aK+bL\right)  E   =q^{2}E\left(  aK+bL\right)  ,\label{ab1}\\
& \left(  a^{-1}\overline{K}+b^{-1}\overline{L}\right)  E
=q^{-2}E\left(
a^{-1}\overline{K}+b^{-1}\overline{L}\right)  ,\\
& \left(  aK+bL\right)  F    =q^{-2}F\left(  aK+bL\right)  ,\\
& \left(  a^{-1}\overline{K}+b^{-1}\overline{L}\right)  F
=q^{2}F\left(
a^{-1}\overline{K}+b^{-1}\overline{L}\right)  ,\\
& EF-FE    =\dfrac{\left(  aK+bL\right)  -\left(  a^{-1}\overline{K}%
+b^{-1}\overline{L}\right)  }{q-q^{-1}} \label{ab5}%
\end{align}
which together with (\ref{k})-(\ref{kk}) and (\ref{l})-(\ref{ll}) determine an algebra we denote by $U_{aK+bL}^{(alg)22}$, the indices $22$ stand for the numbers of
generators in the left (resp., right) hand sides of the relations between the Cartan-like generators ($K$, $L$) and $E$,
$F$. This algebra corresponds to $U_{q}^{w}=\mathfrak{wsl}_{q}\left(
2\right)  $ introduced by Duplij and
Li \cite{dup/li3,dup/li2}. To be more precise, there exists an algebra homomorphism
$\mathfrak{wsl}_{q}\left(
2\right)\rightarrow U_{aK+bL}^{(alg)22}$, which in the notation of \cite{dup/li3}
is given by 
\begin{equation}\label{whom}
K_w\mapsto aK+bL,\qquad \overline{K}_w\mapsto a^{-1}\overline{K}+b^{-1}\overline{L},\qquad
E_w\mapsto E,\qquad F_w\mapsto F.
\end{equation}
As one can see from Lemma \ref{invelt}, together with \eqref{ab1} -- \eqref{ab5},
the image of this homomorphism is a copy of $U_{q}\left(  sl_{2}\right)$,
cf. \cite[Proposition 1]{dup/li3}. 

Next we present an analog of the algebra $U_{q}^{v}=\mathfrak{vsl}
_{q}\left(  2\right)  $ as in \cite{dup/li3}. This is an algebra having the same
generators as $U_{aK+bL}^{(alg)22}$, and being subject to the
relations (together with
(\ref{k}) -- (\ref{kk}) and (\ref{l}) -- (\ref{ll}))%
\begin{align}
& \left(  aK+bL\right)  E\left(
a^{-1}\overline{K}+b^{-1}\overline{L}\right)
  =q^{2}E,\label{abb1}\\
& \left(  a^{-1}\overline{K}+b^{-1}\overline{L}\right)  E\left(
aK+bL\right)
  =q^{-2}E,\\
& \left(  aK+bL\right)  F\left(
a^{-1}\overline{K}+b^{-1}\overline{L}\right)
  =q^{-2}F,\\
& \left(  a^{-1}\overline{K}+b^{-1}\overline{L}\right)  F\left(
aK+bL\right)
  =q^{2}F  ,\label{abb4}
\\
& EF-FE   =\dfrac{\left(  aK+bL\right)  -\left(  a^{-1}\overline{K}
+b^{-1}\overline{L}\right)  }{q-q^{-1}}, \label{abb6}
\end{align}
which we denote $U_{aK+bL}^{(alg)31}$. This algebra corresponds
to the algebra $U_{q}^{v}=\mathfrak{vsl}_{q}\left(  2\right)  $  
  \cite{dup/li3} in the sense that there exists an algebra homomorphism
$\mathfrak{vsl}_{q}\left(
2\right)\rightarrow U_{aK+bL}^{(alg)31}$.
Again, this homomorphism, in the notation of \cite{dup/li3}, is given on
the generators by \eqref{whom}, with the indices $w$ being replaced by $v$.
Another application of Lemma \ref{invelt} allows one to observe that the image of this homomorphism is a copy of $U_{q}\left(  sl_{2}\right)$,
cf. \cite[Proposition 1]{dup/li3}.

 Introduce an extension $\mathbf{\Phi}^{\left(  a,b\right)  }$ of $\mathbf{\Phi}_{\mathcal{H}%
}^{\left(  a,b\right)  }$ to a morphism of $U_{q}\left(  sl_{2}\right)  $  with values in $U_{aK+bL}^{\left(
alg\right)  22}$ and $U_{aK+bL}^{\left(  alg\right)  31}$ as
\begin{equation}
\mathbf{\Phi}^{\left(  a,b\right)  }:\left\{
\begin{array}
[c]{c}%
k\rightarrow aK+bL,\ \ \ \ \ \ k^{-1}\rightarrow a^{-1}\overline{K}
+b^{-1}\overline{L},\\
e\rightarrow E,\ \ \ \ \ f\rightarrow F.
\end{array}
\right.  \label{fab}%
\end{equation}
\begin{proposition}
\label{prop-iso}The algebras $U_{aK+bL}^{\left(  alg\right)  22}$
and $U_{aK+bL}^{\left(  alg\right)  31}$ are isomorphic to
$U_{K+L}^{\left( alg\right)  22}\overset{def}{=}U_{aK+bL}^{\left(
alg\right)  22}|_{a=1,b=1}$ and $U_{K+L}^{\left(  alg\right)
31}\overset{def}{=}U_{aK+bL}^{\left( alg\right)  31}|_{a=1,b=1}$
respectively.
\end{proposition}

\begin{proof}
The desired isomorphism $\mathbf{\Psi}:U_{K+L}^{\left(  alg\right)
22,31}\rightarrow U_{aK+bL}^{\left(  alg\right)  22,31}$ is given by

$ K\rightarrow aK,\;L\rightarrow bL,\;\;\overline{K}\rightarrow
a^{-1}\overline {K},\;\;\overline{L}\rightarrow
b^{-1}\overline{L},\;\;E\rightarrow E,\;\;F\rightarrow F. $
\end{proof}
Therefore, we will not consider the parameters $a$ and $b$ below.

\section{Splitting the relations}

The idempotents $P$ and $Q$ are not central in $U_{K+L}^{\left(
alg\right) 22}$ and $U_{K+L}^{\left(  alg\right)  31}$. By allowing certain
misuse of terminology, we are going
to ''split'' the relations (\ref{ab1}) -- (\ref{ab5}) and
(\ref{abb1}) -- (\ref{abb6}) in such a way that either $P$ and $Q$  become
central
\begin{align}
PE  &  =EP,\ \;QE=EQ,\label{n1}\\
PF  &  =FP,\ \;QF=FQ, \label{n2}%
\end{align}
or satisfy the ``twisting'' conditions
\begin{align}
PE  &  =EQ,\ \;QE=EP,\label{tw1}\\
PF  &  =FQ,\ \;QF=FP. \label{tw2}%
\end{align}
To be more precise, we are about to add the above relations in order to get
the associated quotients of  $U_{K+L}^{\left(alg\right) 22}$ and $U_{K+L}^{\left(  alg\right)  31}$. The ``splitted''\ 22-algebras are given by the
following lists of
relations:
\begin{equation}
\begin{tabular}
[c]{|l|l|}\hline $U_{K,L,norm}^{\left(  alg\right)  22}$ &
$U_{K,L,twist}^{\left(  alg\right) 22}$\\\hline\hline
$K\overline{K}K=K,\quad\overline{K}K\overline{K}=\overline{K}$, &
$K\overline{K}K=K,\quad\overline{K}K\overline{K}=\overline{K}$,\\
$K\overline{K}=\overline{K}K$, & $K\overline{K}=\overline{K}K$,\\
$L\overline{L}L=L,\quad\overline{L}L\overline{L}=\overline{L}$, &
$L\overline{L}L=L,\quad\overline{L}L\overline{L}=\overline{L}$,\\
$L\overline{L}=\overline{L}L$, & $L\overline{L}=\overline{L}L$,\\
$K\overline{K}+L\overline{L}=\mathbf{1}$, & $K\overline{K}+L\overline{L}=\mathbf{1}$,\\
$KE=q^{2}EK,\quad LE=q^{2}EL$, & $KE=q^{2}EL$,\quad $LE=q^{2}EK$,\\
$\overline{K}E=q^{-2}E\overline{K},\quad\overline{L}E=q^{-2}E\overline{L}$, &
$\overline{K}E=q^{-2}E\overline{L},\quad\overline{L}E=q^{-2}E\overline{K}$,\\
$KF=q^{-2}FK$,\quad $LF=q^{-2}FL$, & $KF=q^{-2}FL$,\quad $LF=q^{-2}FK$,\\
$\overline{K}F=q^{2}F\overline{K},\quad\overline{L}F=q^{2}F\overline{L}$, &
$\overline{K}F=q^{2}F\overline{L},\quad\overline{L}F=q^{2}F\overline{K}$,\\
$EF-FE=\dfrac{\left(  K+L\right)  -\left(
\overline{K}+\overline{L}\right) }{q-q^{-1}}$ & $EF-FE=\dfrac{\left(
K+L\right)  -\left(  \overline {K}+\overline{L}\right)
}{q-q^{-1}}$\\\hline
\end{tabular}
\label{22}
\end{equation}
and the ''splitted'' 31-algebras are defined as follows:
\begin{equation}
\begin{tabular}
[c]{|l|l|}\hline $U_{K,L,norm}^{\left(  alg\right)  31}$ &
$U_{K,L,twist}^{\left(  alg\right) 31}$\\\hline\hline
$K\overline{K}K=K,\quad\overline{K}K\overline{K}=\overline{K}$, &
$K\overline{K}K=K,\quad\overline{K}K\overline{K}=\overline{K}$,\\
$K\overline{K}=\overline{K}K$, & $K\overline{K}=\overline{K}K$,\\
$L\overline{L}L=L,\quad\overline{L}L\overline{L}=\overline{L}$, &
$L\overline{L}L=L,\quad\overline{L}L\overline{L}=\overline{L}$,\\
$L\overline{L}=\overline{L}L$, & $L\overline{L}=\overline{L}L$,\\
$K\overline{K}+L\overline{L}=\mathbf{1}$, & $K\overline{K}+L\overline{L}=\mathbf{1}$,\\
$KE\overline{K}=q^{2}EK\overline{K},LE\overline{L}=q^{2}EL\overline{L}$, & $KE\overline{L}=q^{2}EL\overline{L},LE\overline{K}=q^{2}EK\overline{K}$,\\
$\overline{K}EK=q^{-2}EK\overline{K},\overline{L}EL=q^{-2}EL\overline{L}$, & $\overline{K}EL=q^{-2}EL\overline{L},\overline{L}EK=q^{-2}EK\overline{K}$,\\
$KF\overline{K}=q^{-2}FK\overline{K},LF\overline{L}=q^{-2}FL\overline{L}$, &
$KF\overline{L}=q^{-2}FL\overline{L},LF\overline{K}=q^{-2}FK\overline{K}$,\\
$\overline{K}FK=q^{2}FK\overline{K},\overline{L}FL=q^{2}FL\overline{L}$, & $\overline{K}FL=q^{2}FL\overline{L},\overline{L}FK=q^{2}FK\overline{K}$,\\
$K\overline{K}\left(EF-FE\right)=\dfrac{K-\overline{K}}{q-q^{-1}}$, &
$K\overline{K}\left(EF-FE\right)=\dfrac{K-\overline{K}}{q-q^{-1}}$,\\
$L\overline{L}\left(EF-FE\right)=\dfrac{L-\overline{L}}{q-q^{-1}}$ &
$L\overline{L}\left(EF-FE\right)=\dfrac{L-\overline{L}}{q-q^{-1}}$\\\hline
\end{tabular}
\ \label{31}
\end{equation}

Note that $P=K\overline{K}$ and $Q=L\overline{L}$ are not among the generators
used in \eqref{22} and \eqref{31}. The relations which appear in the tables,
  form the (equivalent) translation
in terms of the ''true'' generators of the earlier relations for $U_{K+L}^{\left(alg\right) 22}$ and $U_{K+L}^{\left(  alg\right)  31}$, together with the ''splitting'' relations
\eqref{n1} -- \eqref{tw2}. The procedure of deducing relations in tables from the original ''non-splitted'' relations in most cases reduces to right and/or left multiplication by the idempotents $P$ and $Q$ with subsequent
use of the ''annihilation rules'' \eqref{lk}. Conversely, suppose that \eqref{22}
and \eqref{31} are given. For example, let us start from the relations in
the left column of \eqref{31}. To see that in this case $P$ is central, one
has, using \eqref{lk},
\begin{multline*}
PE=K\overline{K}E(P+Q)=K(\overline{K}EK)\overline{K}+K\overline{K}(EL\overline{L})=\\
=K(q^{-2}EK\overline{K})\overline{K}+K\overline{K}(q^{-2}LE\overline{L})=q^{-2}KE\overline{K}+0=EK\overline{K}=EP.
\end{multline*}
Of course, the similar ideas work also in the rest of verifications.

\begin{proposition}
We have the following isomorphisms: $U_{K,L,norm}^{\left(
alg\right) 22}\cong U_{K,L,norm}^{\left(  alg\right)  31}$, and
$U_{K,L,twist}^{\left( alg\right)  22}\cong U_{K,L,twist}^{\left(
alg\right)  31}$.
\end{proposition}

\begin{proof}
A straightforward computation shows that, in both cases (normal and
twisted), the ideals of relations in question coincide. For
instance, the right multiplication of $KE=q^{2}EK$ by $\overline{K}$
in $U_{K,L,norm}^{\left( alg\right)  22}$ yields
$KE\overline{K}=q^{2}EP$ as in $U_{K,L,norm}^{\left(
alg\right)  31}$. Conversely, starting from the relation $KE\overline{K}%
=q^{2}EP$ in $U_{K,L,norm}^{\left(  alg\right)  31}$ we calculate
$KE=K\left( PE\right)  =K\left(  EP\right)  =\left(
KE\overline{K}\right)  K=\left( q^{2}EP\right)  K=q^{2}EK$ as in
$U_{K,L,norm}^{\left(  alg\right)  22}$. Multiplying the
$EF$-relations in $U_{K,L,norm}^{\left(  alg\right)  22}$,
$U_{K,L,twist}^{\left(  alg\right)  22}$ by $P$ and $Q$ we obtain
the
$EF$-relations of $U_{K,L,norm}^{\left(  alg\right)  31}$, $U_{K,L,twist}%
^{\left(  alg\right)  31}$, and conversely, summing up the last two
$EF$-relations of $U_{K,L,norm}^{\left(  alg\right)  31}$ and using
(\ref{pq1}), we obtain the $EF$-relations of $U_{K,L,norm}^{\left(
alg\right)  22}$. Similar arguments establish the second
isomorphism.
\end{proof}
Therefore, in what follows we consider the algebras
$U_{K,L,norm}^{\left( alg\right)  22}$, $U_{K,L,twist}^{\left(
alg\right)  22}$ (with 22 superscript being discarded) only.

Now we extend the morphism $\mathbf{\Phi}_{\mathcal{H}}$ to that
taking values in the ``splitted'' algebras $U_{K,L,norm}^{\left(
alg\right)  }$ and $U_{K,L,twist}^{\left(  alg\right)  }$ as follows%
\begin{equation}
\mathbf{\Phi}:\left\{
\begin{array}
[c]{c}%
k\rightarrow K+L,\ \ \ \ \ \ k^{-1}\rightarrow\overline{K}+\overline{L},\\
e\rightarrow E,\ \ \ \ f\rightarrow F.
\end{array}
\right.  \label{phi}%
\end{equation}
\begin{proposition}
\label{kerf2}The map $\mathbf{\Phi}$ defined on the generators as
above, admits an extension to a well defined morphism of algebras
from $U_{q}(sl_{2})$ to either $U_{K,L,norm}^{\left(  alg\right)  }$
or $U_{K,L,twist}^{\left( alg\right)  }$, which is an embedding.
\end{proposition}

\begin{proof}
Use an argument similar to that applied in the proof of
\textbf{Proposition \ref{kerf1}}.
\end{proof}

\begin{corollary}
Both algebras $U_{K,L,norm}^{\left(  alg\right)  }$ and $U_{K,L,twist}%
^{\left(  alg\right)  }$ contain $U_{q}\left(  sl_{2}\right)  $ as a
subalgebra.
\end{corollary}

\begin{proof}
Follows from \textbf{Proposition \ref{kerf2}}.
\end{proof}

Note that the Pierce decomposition of $U_{K,L,norm}^{\left(  alg\right)  }$ is%
\begin{equation}
U_{K,L,norm}^{\left(  alg\right)  }=PU_{K,L,norm}^{\left(
alg\right)
}P+QU_{K,L,norm}^{\left(  alg\right)  }Q, \label{upq}%
\end{equation}
which reduces to a direct sum of the two ideals. This leads to

\begin{proposition}
\label{prop-isom}$U_{K,L,norm}^{\left(  alg\right)  }$ is a direct
sum of subalgebras with each summand being isomorphic to
$U_{q}\left(  sl_{2}\right) $.
\end{proposition}

\begin{proof}
The desired isomorphism is given by%
\begin{align}
K  &  \longmapsto k\oplus0,\ \ \overline{K}\longmapsto k^{-1}\oplus
0,\ \ PE\longmapsto e\oplus0,\ PF\longmapsto f\oplus0,\label{plus}\\
L  &  \longmapsto0\oplus k,\ \ \overline{L}\longmapsto0\oplus k^{-1}%
,\ \ QE\longmapsto0\oplus e,\ \ QF\longmapsto0\oplus f,
\end{align}
hence $P\longmapsto\mathbf{1}\oplus0$, $Q$
$\longmapsto0\oplus\mathbf{1}$. This morphism splits as a direct sum
of two morphisms each of the latter being, obviously, an
isomorphism.
\end{proof}

In the ``twisted'' case the Pierce
decomposition%
\begin{equation}
U_{K,L,twist}^{\left(  alg\right)  }=PU_{K,L,twist}^{\left(
alg\right) }P+PU_{K,L,twist}^{\left(  alg\right)
}Q+QU_{K,L,twist}^{\left(  alg\right)
}P+QU_{K,L,twist}^{\left(  alg\right)  }Q, \label{upq1}%
\end{equation}
is nontrivial as all terms are nonzero, i.e. (\ref{upq1}) is not a
direct sum of ideals.

Let us introduce a special automorphism of algebras
$U_{K,L,norm}^{\left( alg\right)  }$ and $U_{K,L,twist}^{\left(
alg\right)  }$, which will be denoted by the same letter
$\mathbf{\Upsilon}$. In either case,
$\mathbf{\Upsilon}$ is given on the generators by%
\begin{equation}
E\mapsto E,~F\mapsto F,~K\mapsto
L,~\overline{K}\mapsto\overline{L},~L\mapsto
K,~\overline{L}\mapsto\overline{K},~\mathbf{1}\mapsto\mathbf{1}, \label{Aut_A}%
\end{equation}
and then extended to an endomorphism of the algebra in question. The
very fact that it becomes this way a well defined linear map and
then its bijectivity is established by observing that
$\mathbf{\Upsilon}$ permutes the list of
generators as well as the list of relations. Note that $\mathbf{\Upsilon}%
^{2}=\mathsf{id}$.

\begin{proposition}\label{prop-norm}
The Poincar\'{e}-Birkhoff-Witt basis of $U_{K,L,norm}^{\left(
alg\right)  }$ is given by the monomials
\begin{align}
&  \left[  \left\{  PK^{i}E^{j}F^{k}\right\}
_{i,j,k\geq0}\cup\left\{
\overline{K}^{i}E^{j}F^{k}\right\}  _{i>0,j,k\geq0}\right] \nonumber\\
&  \cup\left[  \left\{  QL^{i}E^{j}F^{k}\right\}
_{i,j,k\geq0}\cup\left\{ \overline{L}^{i}E^{j}F^{k}\right\}
_{i>0,j,k\geq0}\right]  .
\end{align}
\end{proposition}

\begin{proof}
Since $U_{K,L,norm}^{\left(  alg\right)  }$ is a direct sum of two
copies of $U_{q}(sl_{2})$, the statement immediately follows from
\cite{jantzen}.
\end{proof}

In the case of $U_{K,L,twist}^{\left(  alg\right)  }$ we have the
decomposition into a direct sum of 4 vector subspaces (\ref{upq1}).
We present below a PBW basis which respects this decomposition.

\begin{proposition}
The Poincar\'{e}-Birkhoff-Witt basis of $U_{K,L,twist}^{\left(
alg\right)  }$ is given by the monomials
\begin{align}
&  \left[  \left\{  PK^{i}E^{j}F^{k}\right\}  _{\substack{i,j,k\geq
0\\j+k=even}}\cup\left\{  \overline{K}^{i}E^{j}F^{k}\right\}
_{\substack{i>0,j,k\geq0\\j+k=even}}\right] \nonumber\\
&  \cup\left[  \left\{  PK^{i}E^{j}F^{k}\right\}
_{\substack{i,j,k\geq 0\\j+k=odd}}\cup\left\{
\overline{K}^{i}E^{j}F^{k}\right\}
_{\substack{i>0,j,k\geq0\\j+k=odd}}\right] \nonumber\\
&  \cup\left[  \left\{  QL^{i}E^{j}F^{k}\right\}
_{_{\substack{i,j,k\geq 0\\j+k=odd}}}\cup\left\{
\overline{L}^{i}E^{j}F^{k}\right\}
_{\substack{i>0,j,k\geq0\\j+k=odd}}\right] \nonumber\\
&  \cup\left[  \left\{  QL^{i}E^{j}F^{k}\right\}
_{\substack{i,j,k\geq 0\\j+k=even}}\cup\left\{
\overline{L}^{i}E^{j}F^{k}\right\}
_{\substack{i>0,j,k\geq0\\j+k=even}}\right]  . \label{pbwt}%
\end{align}

\end{proposition}

\begin{proof}
It follows from (\ref{22}) that the linear span of (\ref{pbwt}) is
stable under multiplication by any of the generators $K$,
$\overline{K}$, $L$, $\overline{L}$, $E$, $F$, which implies that
this stability is also valid under multiplication by any element of
$U_{K,L,twist}^{\left(  alg\right)  }$. Since $P$ and $Q$ are among
the basis vectors, this linear span contains $P+Q=\mathbf{1}$, hence
is just the entire algebra. To prove the linear independence of
(\ref{pbwt}) it suffices to prove that every part of this vector
system which is inside a specific Pierce component, is linear
independent. We now stick to the special case of the Pierce
component $P\cdot U_{K,L,twist}^{\left(
alg\right)  }\cdot P$ which is generated by the family of vectors%
\begin{equation}
\left\{  PK^{i}E^{j}F^{k}\right\}  _{\substack{i,j,k\geq0\\j+k=even}%
}\cup\left\{  \overline{K}^{i}E^{j}F^{k}\right\}
_{\substack{i>0,j,k\geq
0\\j+k=even}}, \label{pp}%
\end{equation}
the part of the vector system (\ref{pbwt}) inside the first bracket.
Consider a (finite) linear combination
\begin{equation}
\underset{%
\genfrac{}{}{0pt}{}{i,j,k\geq0}{j+k\quad even}%
}{\sum}\alpha_{ijk}PK^{i}E^{j}F^{k}+\underset{%
\genfrac{}{}{0pt}{}{i>0,~j,k\geq0}{j+k\quad even}%
}{\sum}\beta_{ijk}\overline{K}^{i}E^{j}F^{k} \label{ilk}%
\end{equation}
which is non-trivial (not all $\alpha_{ijk}$ and $\beta_{ijk}$ are
zero). We are about to prove that (\ref{ilk}) is non-zero. For that,
we first use
$\alpha_{ijk}$ and $\beta_{ijk}$ to produce the associated non-trivial linear combination%

\begin{equation}
\underset{%
\genfrac{}{}{0pt}{}{i,j,k\geq0}{j+k\quad even}%
}{\sum}\alpha_{ijk}k^{i}e^{j}f^{k}+\underset{%
\genfrac{}{}{0pt}{}{i>0~,j,k\geq0}{j+k\quad even}%
}{\sum}\beta_{ijk}k^{-i}e^{j}f^{k} \label{alk}%
\end{equation}
in $U_{q}\left(  sl_{2}\right)  $. Since the monomials involved form
a PBW basis in $U_{q}\left(  sl_{2}\right)  $ \cite{jantzen}, the
linear combination (\ref{alk}) is non-zero. Now apply the map
$\mathbf{\Phi}$ (\ref{phi}) to (\ref{alk})
to obtain%
\begin{equation}
\underset{%
\genfrac{}{}{0pt}{}{i,j,k\geq0}{j+k\quad even}%
}{\sum}\alpha_{ijk}\left(  K+L\right)  ^{i}E^{j}F^{k}+\underset{%
\genfrac{}{}{0pt}{}{i>0,~j,k\geq0}{j+k\quad even}%
}{\sum}\beta_{ijk}\left(  \overline{K}+\overline{L}\right)
^{i}E^{j}F^{k}.
\label{elk}%
\end{equation}
As $\mathbf{\Phi}$ is an embedding by \textbf{Proposition}
\textbf{\ref{kerf2}, }we deduce that (\ref{elk}) is non-zero in
$U_{K,L,twist}^{\left(  alg\right)  }$. Observe also that in the
involved monomials $j+k$ is even; it follows that the
projections of (\ref{elk}) to the Pierce components $P\cdot U_{K,L,twist}%
^{\left(  alg\right)  }\cdot Q$ and $Q\cdot U_{K,L,twist}^{\left(
alg\right) }\cdot P$ are both zero. Hence (\ref{elk}) is the sum of
its projections to $P\cdot U_{K,L,twist}^{\left(  alg\right)  }\cdot
P$ and $Q\cdot
U_{K,L,twist}^{\left(  alg\right)  }\cdot Q$, which are just%
\[
\underset{%
\genfrac{}{}{0pt}{}{i,j,k\geq0}{j+k\quad even}%
}{\sum}\alpha_{ijk}PK^{i}E^{j}F^{k}+\underset{%
\genfrac{}{}{0pt}{}{i>0,~j,k\geq0}{j+k\quad even}%
}{\sum}\beta_{ijk}\overline{K}^{i}E^{j}F^{k}%
\]
and%
\[
\underset{%
\genfrac{}{}{0pt}{}{i,j,k\geq0}{j+k\quad even}%
}{\sum}\alpha_{ijk}QL^{i}E^{j}F^{k}+\underset{%
\genfrac{}{}{0pt}{}{i>0,~j,k\geq0}{j+k\quad even}%
}{\sum}\beta_{ijk}\overline{L}^{i}E^{j}F^{k},
\]
respectively. It is easy to see that these are intertwined by the
automorphism $\mathbf{\Upsilon}$ (\ref{Aut_A}), which implies that
these projections are simultaneously zero or non-zero. Of course,
the second assumption is true,
because their sum (\ref{elk}) is non-zero. In particular,%
\[
\underset{%
\genfrac{}{}{0pt}{}{i,j,k\geq0}{j+k\quad even}%
}{\sum}\alpha_{ijk}PK^{i}E^{j}F^{k}+\underset{%
\genfrac{}{}{0pt}{}{i>0,~j,k\geq0}{j+k\quad even}%
}{\sum}\beta_{ijk}\overline{K}^{i}E^{j}F^{k}%
\]
is non-zero, which was to be proved. The proof of linear
independence of all other subsystems of (\ref{pbwt}) (in brackets),
related to other Pierce components, goes in a similar way.
\end{proof}

Let us consider the classical limit $q\rightarrow1$ for
$U_{K,L,norm}^{\left( alg\right)  }$ and $U_{K,L,twist}^{\left(
alg\right)  }$ algebras.

\begin{proposition}
The classical limit of $U_{K,L,norm}^{\left(  alg\right)  }$ is just
a direct sum of two copies of classical limits for $U_{q}\left(
sl_{2}\right)  $ in the sense of \cite{kassel}.
\end{proposition}

\begin{proof}
This follows from \textbf{Proposition \ref{prop-isom}}.
\end{proof}

\section{Hopf algebra structure and von Neumann regular antipode}

To construct a bialgebra we need a counit $\varepsilon$ on
$U_{K+L}$, to be denoted by $\varepsilon$.
Since $P$ and $Q$ are idempotents in $U_{K+L}$, one has
$\varepsilon\left(  P\right)  \left( \varepsilon\left(  P\right)
-1\right)  =0$ and $\varepsilon\left(  Q\right) \left(
\varepsilon\left(  Q\right)  -1\right)  =0$, which implies that
either $\varepsilon\left(  P\right)  =1$, $\varepsilon\left(
Q\right)  =0$ or $\varepsilon\left(  P\right)  =0$,
$\varepsilon\left(  Q\right)  =1$. We assume the first choice. Then
it follows from $L=QL$ that $\varepsilon\left( L\right)
=\varepsilon\left(  QL\right)  =0$. Also it follows from (4) that
$\varepsilon(K+L)=1$, hence $\varepsilon(K)=1$.

Elaborate the embedding $\mathbf{\Phi}$ defined in (\ref{kab}) and
the standard relations (\ref{k0}),(\ref{def}), (\ref{ef}) to
transfer a coproduct onto the
image of $\mathbf{\Phi}$ (\ref{fab}) as follows%
\begin{align}
& \Delta(K+L)    =\left(  K+L\right)  \otimes\left(  K+L\right)  ,
\label{abkl}\\
& \Delta\left(  \overline{K}+\overline{L}\right)     =\left(
\overline {K}+\overline{L}\right)  \otimes\left(
\overline{K}+\overline{L}\right)
,\label{abkl1}\\
& \Delta(E)    =\mathbf{1}\otimes E+E\otimes\left(  K+L\right)  ,\label{abe}\\
& \Delta(F)    =F\otimes\mathbf{1}+\left(
\overline{K}+\overline{L}\right) \otimes
F,\label{abf}\\
& \varepsilon(E)    =\varepsilon(F)=0,\label{eef}\\
& \varepsilon(K+L)    ={1},\label{ekl}\\
& \varepsilon\left(  \overline{K}+\overline{L}\right)     ={1}. \label{ekl1}%
\end{align}

To produce a comultiplication on the above algebras
$U_{K,L,norm}^{\left( alg\right)  }$ and $U_{K,L,twist}^{\left(
alg\right)  }$ determined by (\ref{22}),\medskip\ use
(\ref{abkl})--(\ref{ekl1}) to define a coproduct $\Delta$ first on
$\mathbf{\Phi}\left(  U_{q}^{\left(  alg\right)  }\left(
sl_{2}\right) \right)  $ (via transferring from $U_{q}^{\left(
alg\right) }\left( sl_{2}\right)  $) and then extend it to the
entire algebras $U_{K,L,norm}^{\left(  alg\right)  }$ and
$U_{K,L,twist}^{\left( alg\right)
}$ as follows.%

\begin{equation}%
\begin{tabular}
[c]{|l|l|}\hline $U_{K,L,norm}^{\left(  coalg\right)  }$ &
$U_{K,L,twist}^{\left( coalg\right)  }$\\\hline\hline
$\Delta(K)=K\otimes K,$ & $\Delta(K)=K\otimes K+L\otimes L,$\\
$\Delta(\overline{K})=\overline{K}\otimes\overline{K},$ &
$\Delta(\overline
{K})=\overline{K}\otimes\overline{K}+\overline{L}\otimes\overline{L},$\\
$\Delta(L)=L\otimes L+L\otimes K+K\otimes L,$ & $\Delta(L)=L\otimes
K+K\otimes
L,$\\
$\Delta(\overline{L})=\overline{L}\otimes\overline{L}+\overline{L}%
\otimes\overline{K}+\overline{K}\otimes\overline{L},$ &
$\Delta(\overline
{L})=\overline{L}\otimes\overline{K}+\overline{K}\otimes\overline{L}$\\
$\Delta(E)=\mathbf{1}\otimes E+E\otimes\left(  K+L\right)  ,$ &
$\Delta(E)=\mathbf{1}\otimes
E+E\otimes\left(  K+L\right)  ,$\\
$\Delta(F)=F\otimes\mathbf{1}+\left(
\overline{K}+\overline{L}\right) \otimes F,$ &
$\Delta(F)=F\otimes\mathbf{1}+\left(  \overline{K}+\overline{L}\right)  \otimes F,$\\
$\varepsilon(E)=\varepsilon(F)=0,$ & $\varepsilon(E)=\varepsilon(F)=0,$\\
$\varepsilon(K)=1,\varepsilon(\overline{K})=1,$ & $\varepsilon
(K)=1,\varepsilon(\overline{K})=1,$\\
$\varepsilon(L)=\varepsilon(\overline{L})=0.$ &
$\varepsilon(L)=\varepsilon (\overline{L})=0.$\\\hline
\end{tabular}\label{coalg}
\end{equation}
The convolution on the bialgebras $U_{K,L,norm}^{\left(
bialg\right)  }$and
$U_{K,L,twist}^{\left(  bialg\right)  }$ produced this way is defined by%
\begin{equation}
\left(  \mathsf{A}\star\mathsf{B}\right)  \equiv\mu\left(  \mathsf{A}%
\otimes\mathsf{B}\right)  \Delta, \label{ab}%
\end{equation}
where $\mathsf{A}$\textsf{,}$\mathsf{B}$ are linear endomorphisms of
the underlying vector space.

Let us first consider the bialgebra $U_{K,L,norm}^{\left(
bialg\right)  }$ from viewpoint of Hopf algebra structure.

\begin{proposition}
\label{prop-nos}The bialgebra $U_{K,L,norm}^{\left(  bialg\right)
}$ has no conventional antipode $\mathsf{S}$ satisfying the standard
Hopf algebra axiom
\begin{equation}
\mathsf{S}\star\mathsf{id}=\mathsf{id}\star\mathsf{S}=\eta
\circ\varepsilon. \label{s}%
\end{equation}

\end{proposition}

\begin{proof}
Since $\varepsilon\left(  P\right)  =1$ and $\Delta(P)=P\otimes P$
we have
from (\ref{ab})%
\begin{equation}
\left(  \mathsf{S}\star\mathsf{id}\right)  \left(  P\right)  =\mathsf{S}%
\left(  P\right)  P=\left(  \mathsf{id}\star\mathsf{S}\right) \left(
P\right)  =P\mathsf{S}\left(  P\right)
=\mathbf{1}\cdot\varepsilon\left( P\right)
=\mathbf{1}, \label{s1}%
\end{equation}
which is impossible since $P$ is not invertible.
\end{proof}
Let us introduce an antimorphism $\mathsf{T}$ of
$U_{K,L,norm}^{\left(
bialg\right)  }$ as follows%
\begin{align}
\mathsf{T}\left(  K\right)   &  =\overline{K},\;\mathsf{T}\left(
\overline {K}\right)  =K,\;\mathsf{T}\left(  L\right)
=\overline{L},\;\mathsf{T}\left(
\overline{L}\right)  =L,\;\label{t1}\\
\mathsf{T}\left(  E\right)   &  =-E\left(
\overline{K}+\overline{L}\right)
,\ \ \mathsf{T}\left(  F\right)  =-\left(  K+L\right)  F. \label{t2}%
\end{align}
For $U_{K,L,norm}^{\left(  bialg\right)  }$ we observe that%
\begin{align}
\left(  \mathsf{T}\star\mathsf{id}\right)  \left(  K\right)   &
=\left( \mathsf{id}\star\mathsf{T}\right)  \left(  K\right) =\left(
\mathsf{T}\star\mathsf{id}\right)  \left( \overline{K}\right)
=\left( \mathsf{id}\star\mathsf{T}\right) \left( \overline{K}\right)
=P,\label{idn1}\\
\left(  \mathsf{T}\star\mathsf{id}\right)  \left(  L\right)   &
=\left( \mathsf{id}\star\mathsf{T}\right)  \left(  L\right) =\left(
\mathsf{T}\star\mathsf{id}\right)  \left( \overline{L}\right)
=\left( \mathsf{id}\star\mathsf{T}\right) \left( \overline{L}\right)
=Q,\label{idn2}\\
 \left(
\mathsf{T}\star\mathsf{id}\right)  \left(  E\right)  & = \left(
\mathsf{id}\star\mathsf{T}\right)  \left(  E\right)=\left(
\mathsf{T}\star\mathsf{id}\right)  \left(  F\right) =\left(
\mathsf{id}\star\mathsf{T}\right)  \left(  F\right)  =0.\label{idn3}%
\end{align}

\begin{proposition}
The antimorphism $\mathsf{T}$ of \ $U_{K,L,norm}^{\left(
bialg\right) }$ is
von Neumann regular%
\begin{equation}
\mathsf{id}\star\mathsf{T}\star\mathsf{id}=\mathsf{id},\ \ \ \ \ \mathsf{T}%
\star\mathsf{id}\star\mathsf{T}=\mathsf{T}.\label{tt}%
\end{equation}

\end{proposition}

\begin{proof}
First observe that, since a convolution of linear maps is again a
linear map, it suffices to verify (\ref{tt}) separately on the
direct summands $PU_{K,L,norm}^{\left(  bialg\right)  }$ and
$QU_{K,L,norm}^{\left( bialg\right)  }$, associated to the central
idempotents $P$ and $Q$, respectively. We start with
$PU_{K,L,norm}^{\left(  bialg\right)  }$, which is a sub-bialgebra.
Denote by $\varphi_{P}:PU_{K,L,norm}^{\left(  bialg\right)
}\rightarrow U_{q}\left(  sl_{2}\right)  $ the isomorphism (\ref{plus}).
Earlier it was introduced as an isomorphism of algebras (hence it
intertwines the products, $\varphi_{P}\circ\mu\circ\left(
\varphi_{P}^{-1}\otimes\varphi _{P}^{-1}\right)  =\mu_0=\mu_{U_{q}\left(
sl_{2}\right)  }$), but now it follows from (\ref{coalg}) and
$\Delta(P)=P\otimes P$ that $\varphi_{P}$ also intertwines the
comultiplication (\ref{k0})-(\ref{def}) of $U_{q}\left(
sl_{2}\right)  $ and the restriction of the comultiplication
$\Delta$ of $U_{K,L,norm}^{\left( bialg\right)  }$ onto
$PU_{K,L,norm}^{\left(  bialg\right)  }$, that is, $\left(
\varphi_{P}\otimes\varphi_{P}\right)  \circ\Delta\circ\varphi
_{P}^{-1}=\Delta_{0}$.

It follows that, given any two endomorphisms of the underlying
vector space of $U_{K,L,norm}^{\left(  bialg\right)  }$which leave
$PU_{K,L,norm}^{\left( bialg\right)  }$ invariant, then
$\varphi_{P}$ sends the convolution of them (restricted to
$PU_{K,L,norm}^{\left(  bialg\right)  }$) to the convolution of the
transferred maps on $U_{q}\left(  sl_{2}\right)  $.

An obvious verification shows that both $\mathsf{id}$ and
$\mathsf{T}$ leave $PU_{K,L,norm}^{\left(  bialg\right)  }$
invariant, and then a computation
shows that so do $\mathsf{id}\star\mathsf{T}$ and $\mathsf{T}\star\mathsf{id}%
$. Specifically, one has%
\[
\left(  \mathsf{id}\star\mathsf{T}\right)  \left(  PX\right)  \mathsf{=}%
\left(  \mathsf{T}\star\mathsf{id}\right)  \left(  PX\right)
=\varepsilon_0\left(  \varphi_{P}\left(  PX\right)  \right)  P
\]
for any $X\in U_{K,L,norm}^{\left(  bialg\right)  }$. This means
that
$\varphi_{P}$ establishes the equivalence of (\ref{tt}) on $PU_{K,L,norm}%
^{\left(  bialg\right)  }$ and the von Neumann regularity conditions
for the transfer of $\mathsf{T}$ via $\varphi_{P}$ on $U_{q}\left(
sl_{2}\right)  $. An easy verification shows that this transfer
is just $\mathsf{S}$, the antipode of $U_{q}\left(  sl_{2}\right)
$. It is well known that $\mathsf{S}$ is also von Neumann regular,
which finishes the proof of (\ref{tt}) restricted to
$PU_{K,L,norm}^{\left(  bialg\right)  }$.

On can readily replace in the above argument $\varphi_{P}$ by the
isomorphism $\mathbf{\Phi}^{-1}:\mathbf{\Phi}\left(  U_{q}\left(
sl_{2}\right) \right)  \rightarrow U_{q}\left(  sl_{2}\right)  $,
with $\mathbf{\Phi}$ being the embedding (\ref{phi}). This way we
obtain (\ref{tt}) restricted to $\mathbf{\Phi}\left(  U_{q}\left(
sl_{2}\right) \right)  $. However, this argument is inapplicable to
$QU_{K,L,norm}^{\left( bialg\right)  }$, as the latter fails to be a
sub-coalgebra.

Now observe that the projection of $\mathbf{\Phi}\left(  U_{q}\left(
sl_{2}\right) \right)  $ to the direct summand
$QU_{K,L,norm}^{\left(  bialg\right)  }$ is just
$QU_{K,L,norm}^{\left(  bialg\right)  }$. This is because the PBW
basis $\left\{  k^{i}e^{j}f^{k}\right\}  _{j,k\geq0}$ of
$U_{q}\left(
sl_{2}\right)  $ transferred by $\mathbf{\Phi}$ is just%
\[
\left\{  \left(  K+L\right)  ^{i}E^{j}F^{k}\right\}
_{i,j,k\geq0}\cup\left\{ \left(  \overline{K}+\overline{L}\right)
^{i}E^{j}F^{k}\right\} _{i>0,j,k\geq0}.
\]
These vectors project to $QU_{K,L,norm}^{\left(  bialg\right)  }$ as%
\[
\left\{  QL^{i}E^{j}F^{k}\right\}  _{i,j,k\geq0}\cup\left\{  \overline{L}%
^{i}E^{j}F^{k}\right\}  _{i>0,j,k\geq0},
\]
which form a basis in $QU_{K,L,norm}^{\left(  bialg\right)  }$ by
\textbf{Proposition \ref{prop-norm}}. Thus, given any $X\in U_{K,L,norm}%
^{\left(  bialg\right)  }$, one can find $x\in U_{q}\left(
sl_{2}\right)  $
such that $QX=Q\mathbf{\Phi}\left(  x\right)  $. In view of this, one has%
\begin{align*}
&  \left(  \mathsf{id}\star\mathsf{T}\star\mathsf{id}\right)  \left(
QX\right)  =
  \left(  \mathsf{id}\star\mathsf{T}\star\mathsf{id}\right)  \left(
\left( \mathbf{1}-P\right)  \mathbf{\Phi}\left(  x\right)  \right)
\\&=
  \left(  \mathsf{id}\star\mathsf{T}\star\mathsf{id}\right)  \left(
\mathbf{\Phi}\left(  x\right)  \right)  -\left(  \mathsf{id}\star\mathsf{T}%
\star\mathsf{id}\right)  \left(  P\mathbf{\Phi}\left(  x\right)  \right)  \\
& = \mathbf{\Phi}\left(  x\right)  -P\mathbf{\Phi}\left(  x\right)
=\left( \mathbf{1}-P\right) \mathbf{\Phi}\left(  x\right)
=Q\mathbf{\Phi}\left( x\right)  =
  QX
\end{align*}
due to the above observations. Certainly, a similar computation is
applicable to the second part of (\ref{tt}), which completes its
verification on $QU_{K,L,norm}^{\left(  bialg\right)  }$, hence on
$U_{K,L,norm}^{\left( bialg\right)  }$.
\end{proof}

\begin{definition}
We call the antimorphism $\mathsf{T}$ with property (\ref{tt}) a von
Neumann regular antipode.
\end{definition}

\begin{definition}
We call a bialgebra with a von Neumann regular antipode a von
Neumann-Hopf algebra.
\end{definition}

\begin{remark}
The standard Drinfeld-Jimbo algebra $U_{q}\left(  sl_{2}\right)  $
(which is a domain \cite{jantzen}) admits no embedding of
$U_{K,L,norm}^{\left( bialg\right)  }$, because the latter contain
zero divisors (e.g. (\ref{pq1})).
\end{remark}

Let us consider a possibility to produce a Hopf algebra structure on
$U_{K,L,twist}^{\left(  bialg\right)  }$. First we observe that the
argument of the proof of \textbf{Proposition \ref{prop-nos}} does
not work in this case. Indeed, an application of (\ref{s}) to $P$
yields, instead of
(\ref{s1}), the following relation%
\begin{equation}
\mathsf{S}\left(  P\right)  P+\mathsf{S}\left(  Q\right)  Q=\mathbf{1}, \label{spq}%
\end{equation}
which does not contradict to noninvertibility of $P$ and $Q$ as in
the context
of (\ref{s1}). Introduce an antimorphism $\mathsf{S}$ of $U_{K,L,twist}%
^{\left(  bialg\right)  }$ by the same formulas as (\ref{t1})--(\ref{t2})%
\begin{align}
\mathsf{S}\left(  K\right)   &  =\overline{K},\mathsf{S}\left(
\overline {K}\right)  =K,\mathsf{S}\left(  L\right)
=\overline{L},\mathsf{S}\left(
\overline{L}\right)  =L,\\
\mathsf{S}\left(  E\right)   &  =-E\left(
\overline{K}+\overline{L}\right) ,\ \ \mathsf{S}\left(  F\right)
=-\left(  K+L\right)  F.
\end{align}

We have for $U_{K,L,twist}^{\left(  bialg\right)  }$%
\begin{align}
\left(  \mathsf{id}\star\mathsf{S}\right)  \left(  K\right)   &
=\left( \mathsf{S}\star\mathsf{id}\right)  \left(  K\right) =\left(
\mathsf{S}\star\mathsf{id}\right)  \left( \overline{K}\right)
=\left( \mathsf{id}\star\mathsf{S}\right) \left( \overline{K}\right)
=\mathbf{1},\label{idt1}\\
\left(  \mathsf{id}\star\mathsf{S}\right)  \left(  L\right)   &
=\left( \mathsf{S}\star\mathsf{id}\right)  \left(  L\right) =\left(
\mathsf{S}\star\mathsf{id}\right)  \left( \overline{L}\right)
=\left( \mathsf{id}\star\mathsf{S}\right) \left( \overline{L}\right)
=0,\label{idt2}\\
\left(  \mathsf{id}\star\mathsf{S}\right)  \left(  E\right)   &
=\left( \mathsf{S}\star\mathsf{id}\right)  \left(  E\right) =\left(
\mathsf{S}\star\mathsf{id}\right)  \left(  F\right) =\left(
\mathsf{id}\star\mathsf{S}\right)  \left(  F\right)  =0.\label{idt3}%
\end{align}

The proof of the following statement is basically due to
\cite[p.35]{jantzen}.

\begin{proposition}
The relations $\left(  \mathsf{id}\star\mathsf{S}\right)  \left(
X\right) =\left(  \mathsf{S}\star\mathsf{id}\right)  \left( X\right)
=\varepsilon\left(  X\right)  \cdot\mathbf{1}$ are valid for any $X\in U_{K,L,twist}%
^{\left(  bialg\right)  }$.
\end{proposition}

\begin{proof}Note that $X\mapsto\varepsilon{X}\mathbf{1}$ is a morphism of algebras. Hence,
in view of an obvious induction argument, it suffices to verify that
$\left( \mathsf{id}\star\mathsf{S}\right)  \left(  XY\right) =\left(
\mathsf{id}\star\mathsf{S}\right)  \left(  X\right) \cdot\left(
\mathsf{id}\star\mathsf{S}\right)  \left(  Y\right) $ and $\left(
\mathsf{S}\star\mathsf{id}\right)  \left( XY\right)  =\left(
\mathsf{S}\star\mathsf{id}\right)  \left( X\right)  \cdot\left(
\mathsf{S}\star\mathsf{id}\right)  \left( Y\right)  $, with $X$
being one of the generators $K,\overline{K},L,\overline{L},E,F$ and
$Y$ arbitrary. We use the Sweedler notation $\Delta\left(  X\right)
=\sum_{i}X_{i}^{\prime
}\otimes X_{i}^{\prime\prime}$ \cite{sweedler} to get%
\[
\left(  \mathsf{S}\star\mathsf{id}\right)  \left(  XY\right)  =\sum
_{ij}\mathsf{S}\left(  Y_{j}^{\prime}\right)  \mathsf{S}\left(
X_{i}^{\prime }\right)  X_{i}^{\prime\prime}Y_{j}^{\prime\prime}.
\]
It follows from (\ref{idt1})--(\ref{idt3}) that
$\sum_{i}\mathsf{S}\left( X_{i}^{\prime}\right)
X_{i}^{\prime\prime}$ is a scalar multiple of $\mathbf{1}$,
hence is central in $U_{K,L,twist}^{\left(  bialg\right)  }$, and we obtain%
\begin{align*}
&\left(  \mathsf{S}\star\mathsf{id}\right)  \left(  XY\right)
=\sum_{ij}\mathsf{S}\left(  X_{i}^{\prime}\right)
X_{i}^{\prime\prime
}\mathsf{S}\left(  Y_{j}^{\prime}\right)  Y_{j}^{\prime\prime}\\
&  =\left(  \sum_{i}\mathsf{S}\left(  X_{i}^{\prime}\right)  X_{i}%
^{\prime\prime}\right)  \left(  \sum_{j}\mathsf{S}\left(
Y_{j}^{\prime }\right)  Y_{j}^{\prime\prime}\right)
  =\left(  \mathsf{S}\star\mathsf{id}\right)  \left(  X\right)
\cdot\left(  \mathsf{S}\star\mathsf{id}\right)  \left(  Y\right) .
\end{align*}
Of course, a similar argument goes also for $\left(  \mathsf{id}%
\star\mathsf{S}\right)  $.
\end{proof}

Thus, we have the following

\begin{theorem}
1) $U_{K,L}^{\left(  Hopf\right)  }\overset{def}{=}\left(  U_{K,L,twist}%
^{\left(  bialg\right)  },\mathsf{S}\right)  $ is a Hopf algebra;

2) $U_{K,L}^{\left(  vN-Hopf\right)  }\overset{def}{=}\left(  U_{K,L,norm}%
^{\left(  bialg\right)  },\mathsf{T}\right)  $ is a von Neumann-Hopf
algebra.
\end{theorem}

\section{Structure of $R$-matrix and the Pierce decomposition}

Let us consider a version of universal $R$-matrix for
$U_{K,L}^{\left( vN-Hopf\right)  }$ and $U_{K,L}^{\left(
Hopf\right)  }$. In order to avoid considerations related to formal
series (the general context of $R$-matrices), we turn to
quasi-cocommutative bialgebras \cite{kassel}. Such bialgebras
generate $R$-matrices of some simpler shape admitting (under some
additional assumptions) an explicit formula to be described below.

\begin{definition}
\label{def-r}A bialgebra $U^{\left(  bialg\right)  }=\left(  \mathbb{C}%
,B,\mu,\eta,\Delta,\varepsilon\right)  $ is called
quasi-cocommutative, if there exists an invertible element $R\in
U^{\left(  bialg\right)  }\otimes
U^{\left(  bialg\right)  }$, called a universal $R$-matrix, such that%
\begin{equation}
\Delta^{cop}\left(  b\right)  =R\Delta\left(  b\right)  R^{-1},\ \ \
\ \forall
b\in U^{\left(  bialg\right)  }, \label{op}%
\end{equation}
where $\Delta^{cop}$ is the opposite comultiplication in $U^{\left(
bialg\right)  }$.
\end{definition}

The $R$-matrix of a braided bialgebra $U^{\left(  bialg\right)  }$
is subject
to%
\begin{equation}
(\Delta\otimes\operatorname{id})(R)=R_{13}R_{23},\ \ (\operatorname{id}%
\otimes\Delta)(R)=R_{13}R_{12}, \label{rr}%
\end{equation}
where for $R=\sum_{i}s_{i}\otimes t_{i}$ one has
$R_{12}=\sum_{i}s_{i}\otimes t_{i}\otimes\mathbf{1}$, etc.
\cite{dri2}. From now on we assume that $q^{n}=1$, which is a
distinct case to the above context.

Consider the two-sided ideal $I_{sl_{2}}$ in $U_{q}^{\left(
alg\right) }\left(  sl_{2}\right)  $ generated by $\left\{
k^{n}-\mathbf{1},e^{n},f^{n}\right\} $, together with the associated
quotient algebra $\widehat{U}_{q}^{\left( alg\right)  }\left(
sl_{2}\right) =U_{q}^{\left(  alg\right)  }\left( sl_{2}\right)
\diagup I_{sl_{2}}$.

\begin{theorem}
[{\cite[p.230]{kassel}}]\label{theor-kas}The universal $R$-matrix of
$\widehat{U}_{q}^{\left(  alg\right)  }\left(  sl_{2}\right)  $ is%
\begin{align}
\widehat{R}  &  =\sum_{0\leq i,j,m\leq n-1}A_{m}^{ij}\left(
q\right)  \cdot
e^{m}k^{i}\otimes f^{m}k^{j},\label{ra}\\
A_{m}^{ij}\left(  q\right)   &  =\frac{1}{n}\frac{(q-q^{-1})^{m}}%
{[m]!}q^{\tfrac{m(m-1)}{2}+2m(i-j)-2ij}, \label{aq}%
\end{align}
where $[m]!=\left[  1\right]  \left[  2\right]  \ldots\left[
m\right]  $,
$\left[  m\right]  =\left(  q^{m}-q^{-m}\right)  \diagup\left(  q-q^{-1}%
\right)  .$\medskip
\end{theorem}

Now we use (\ref{phi}) to obtain an analog of this theorem for $U_{K,L}%
^{\left(  Hopf\right)  }$. In a similar way we consider the quotient
algebra $\widehat{U}_{K+L}^{\left(  Hopf\right)  }=U_{K,L}^{\left(
Hopf\right) }\diagup I_{K+L}^{\left(  Hopf\right)  }$, where the
two-sided ideal
$I_{K+L}^{\left(  Hopf\right)  }$ is generated by $\left\{  K^{n}%
+L^{n}-\mathbf{1},E^{n},F^{n}\right\}  $.

\begin{theorem}
\label{theo-hopf}The universal $R$-matrix of
$\widehat{U}_{K,L}^{\left(
Hopf\right)  }$ is given by%
\begin{equation}
\widehat{R}_{K+L}^{\left(  Hopf\right)  }=\sum_{0\leq i,j,m\leq n-1}A_{m}%
^{ij}\left(  q\right)  \cdot E^{m}\left(  K^{i}+L^{i}\right)
\otimes F^{m}\left(  K^{j}+L^{j}\right)  .
\end{equation}

\end{theorem}

\begin{proof}
In view of the morphism
$\widehat{\mathbf{\Phi}}:\widehat{U}_{q}^{\left( alg\right) }\left(
sl_{2}\right) \rightarrow\widehat{U}_{K+L}^{\left(  Hopf\right)  }$
induced by (\ref{phi}) and \textbf{Theorem \ref{theor-kas}}, it
suffices (due to invertibility of $R$) to verify the relation
$\Delta^{cop}\left( b\right) \widehat{R}_{K+L}^{\left(  Hopf\right)
}=\widehat{R}_{K+L}^{\left( Hopf\right)  }\Delta\left(  b\right)  $
for $b=K,\overline{K}$, because $\Delta$ and $\Delta^{cop}$ are
morphisms of algebras. This claim reduces to the verification of
\begin{align}
&  \left(  K\otimes K+L\otimes L\right)  \left(  E^{m}\left(  K^{i}%
+L^{i}\right)  \otimes F^{m}\left(  K^{j}+L^{j}\right)  \right) \nonumber\\
&  =\left(  E^{m}\left(  K^{i}+L^{i}\right)  \otimes F^{m}\left(  K^{j}%
+L^{j}\right)  \right)  \left(  K\otimes K+L\otimes L\right)  ,
\end{align}
and%
\begin{align}
&  \left(
\overline{K}\otimes\overline{K}+\overline{L}\otimes\overline
{L}\right)  \left(  E^{m}\left(
\overline{K}^{i}+\overline{L}^{i}\right) \otimes F^{m}\left(
\overline{K}^{j}+\overline{L}^{j}\right)  \right)
\nonumber\\
&  =\left(  E^{m}\left(  \overline{K}^{i}+\overline{L}^{i}\right)
\otimes F^{m}\left(  \overline{K}^{j}+\overline{L}^{j}\right)
\right)  \left(
\overline{K}\otimes\overline{K}+\overline{L}\otimes\overline{L}\right)
,
\end{align}
using (\ref{22}). The relations (\ref{rr}) are transferred by
$\widehat{\mathbf{\Phi}}$ into our picture, because
$\widehat{R}_{K+L}^{\left(  Hopf\right)  }$ is inside of the tensor
square of the image of $\widehat{\mathbf{\Phi}}$.
\end{proof}

Turn to writing down an explicit form for the universal $R$-matrix
in the case of $U_{K,L}^{\left(  vN-Hopf\right)  }$. Again we
consider the quotient algebra $\widehat{U}_{K+L}^{\left(
vN-Hopf\right)  }=U_{K,L}^{\left( vN-Hopf\right)  }\diagup
I_{K+L}^{\left(  vN-Hopf\right)  }$, where the two-sided ideal
$I_{K,L}^{\left(  vN-Hopf\right)  }$ is generated by $\left\{
K^{n}+L^{n}-\mathbf{1},E^{n},F^{n}\right\}  $.

\begin{theorem}
The universal $R$-matrix of $\widehat{U}_{K+L}^{\left(
vN-Hopf\right)  }$ is
given by%
\begin{equation}
\widehat{R}_{K+L}^{\left(  vN-Hopf\right)  }=\sum_{0\leq i,j,m\leq n-1}%
A_{m}^{ij}\left(  q\right)  \cdot E^{m}\left(  K^{i}+L^{i}\right)
\otimes F^{m}\left(  K^{j}+L^{j}\right)  .
\end{equation}
\end{theorem}

\begin{proof}
Is the same as that of \textbf{Theorem \ref{theo-hopf}}.
\end{proof}

\begin{remark}
In view of \textbf{Theorem \ref{theor-kas}} the $R$-matrices we have
introduced satisfy the Yang-Baxter equation by our construction.
\end{remark}

Note that $\widehat{R}_{K+L}^{\left(  vN-Hopf\right)  }$ is not
submitted to the direct sum decomposition (\ref{upq}). Now we
present another notion of $R$-matrix which  respects (\ref{upq}),
but differs from that described in \textbf{Definition \ref{def-r}}
in the sense of being noninvertible.

\begin{definition}
A bialgebra $\widetilde{U}^{\left(  bialg\right)  }=\left(  \mathbb{C}%
,B,\mu,\eta,\Delta,\varepsilon\right)  $ is called
near-quasi-cocommutative, if there exists an element
$\widetilde{R}\in\widetilde{U}^{\left( bialg\right)
}\otimes\widetilde{U}^{\left(  bialg\right)  }$, called a
universal near-$R$-matrix, such that%
\begin{equation}
\Delta^{cop}\left(  b\right)
\widetilde{R}=\widetilde{R}\Delta\left( b\right)  ,\ \ \ \ \forall
b\in\widetilde{U}^{\left(  bialg\right)  },
\end{equation}
where $\Delta^{cop}$ is the opposite comultiplication in $\widetilde
{U}^{\left(  bialg\right)  }$ and an element $\widetilde{R}^{\dagger}%
\in\widetilde{U}^{\left(  bialg\right)
}\otimes\widetilde{U}^{\left(
bialg\right)  }$ is such that%
\begin{equation}
\widetilde{R}\widetilde{R}^{\dagger}\widetilde{R}=\widetilde{R},\ \
\widetilde
{R}^{\dagger}\widetilde{R}\widetilde{R}^{\dagger}=\widetilde{R}^{\dagger},
\label{rrr}%
\end{equation}
and $\widetilde{R}^{\dagger}$ can be named the Moore-Penrose inverse
for a near-$R$-matrix \cite{nashed,rao/mit}.
\end{definition}

A near-quasi-cocommutative bialgebra $\widetilde{U}^{\left(
bialg\right)  }$ is braided, if its near-$R$-matrix satisfies
(\ref{rr}).

Consider the quotient algebra $\widehat{U}_{K,L}^{\left(
vN-Hopf\right) }=U_{K,L}^{\left(  vN-Hopf\right)  }\diagup
I_{K,L}^{\left(  vN-Hopf\right) }$, where the two-sided ideal
$I_{K,L}^{\left(  vN-Hopf\right)  }$ is generated by $\left\{
K^{n}-P,L^{n}-Q,E^{n},F^{n}\right\}  $.

\begin{theorem}
The universal $R$-matrix of $\widehat{U}_{K,L}^{\left(
vN-Hopf\right)  }$ is
given by the sum%
\begin{equation}
\widehat{R}_{K,L}^{\left(  vN-Hopf\right)
}=\widehat{R}_{PP}^{\left( vN-Hopf\right)
}+\widehat{R}_{QQ}^{\left(  vN-Hopf\right)  },
\end{equation}
where
\begin{align}
\widehat{R}_{PP}^{\left(  vN-Hopf\right)  }  &  =\sum_{0\leq
i,j,m\leq
n-1}A_{m}^{ij}\left(  q\right)  \cdot E^{m}K^{i}\otimes F^{m}K^{j},\\
\widehat{R}_{QQ}^{\left(  vN-Hopf\right)  }  &  =\sum_{0\leq
i,j,m\leq n-1}A_{m}^{ij}\left(  q\right)  \cdot E^{m}L^{i}\otimes
F^{m}L^{j}.
\end{align}

\end{theorem}

\begin{remark}
The universal near-$R$-matrix $\widehat{R}_{K,L}^{\left(
vN-Hopf\right)  }$
can be presented in the form%
\begin{equation}
\widehat{R}_{K,L}^{\left(  vN-Hopf\right)  }=\left(  P\otimes
P\right) \widehat{R}_{PP}^{\left(  vN-Hopf\right)  }+\left(
Q\otimes Q\right)
\widehat{R}_{QQ}^{\left(  vN-Hopf\right)  }. \label{rpq}%
\end{equation}

\end{remark}

\begin{proof}
Recall that $U_{K,L}^{\left(  vN-Hopf\right)  }$ admits the direct
sum decomposition (\ref{upq}) with each summand being isomorphic to
$U_{q}\left( sl_{2}\right)  $. After dividing out by the ideal
$I_{K,L}^{\left( vN-Hopf\right)  }$ we get
\begin{align}
\widehat{U}_{K,L}^{\left(  vN-Hopf\right)  }  &  =PU_{K,L}^{\left(
vN-Hopf\right)  }P\diagup\left\{  I_{K,L}^{\left(  vN-Hopf\right)
}\cap
PU_{K,L}^{\left(  vN-Hopf\right)  }P\right\} \nonumber\\
&  +QU_{K,L}^{\left(  vN-Hopf\right)  }Q\diagup\left\{
I_{K,L}^{\left( vN-Hopf\right)  }\cap QU_{K,L}^{\left(
vN-Hopf\right)  }Q\right\}  .
\label{us}%
\end{align}
Each of the summands of the right hand side of (\ref{us}) is clearly
isomorphic to $\widehat{U}_{q}^{\left(  alg\right)  }\left(
sl_{2}\right)  $, and the isomorphisms in question take
$\mathbf{1}\in\widehat{U}_{q}^{\left( alg\right)  }\left(
sl_{2}\right)  $ to $P$ and $Q$ respectively. Now it follows from
\textbf{Theorem \ref{theor-kas}}, that each of the terms of
(\ref{rpq}) satisfies the conditions of \textbf{Definition
\ref{def-r}} and (\ref{rr}), hence so does their sum
$\widehat{R}_{K,L}^{\left( vN-Hopf\right)  }$. Also it follows from
\textbf{Theorem \ref{theor-kas}}, that there exist
$\widehat{R}_{PP}^{\left(  vN-Hopf\right)  \dagger}$,
$\widehat{R}_{QQ}^{\left(  vN-Hopf\right)  \dagger}$ $\in\widehat{U}%
_{K,L}^{\left(  vN-Hopf\right)  }\otimes\widehat{U}_{K,L}^{\left(
vN-Hopf\right)  }$ such that%
\begin{align}
\widehat{R}_{PP}^{\left(  vN-Hopf\right)  }\widehat{R}_{PP}^{\left(
vN-Hopf\right)  \dagger}  &  =\widehat{R}_{PP}^{\left(
vN-Hopf\right)
\dagger}\widehat{R}_{PP}^{\left(  vN-Hopf\right)  }=P\otimes P,\\
\widehat{R}_{QQ}^{\left(  vN-Hopf\right)  }\widehat{R}_{QQ}^{\left(
vN-Hopf\right)  \dagger}  &  =\widehat{R}_{QQ}^{\left(
vN-Hopf\right) \dagger}\widehat{R}_{QQ}^{\left(  vN-Hopf\right)
}=Q\otimes Q,
\end{align}
hence the von Neumann regularity (\ref{rrr}) is valid for
\begin{equation}
\widehat{R}^{\left(  vN-Hopf\right)  }=\widehat{R}_{PP}^{\left(
vN-Hopf\right)  }+\widehat{R}_{QQ}^{\left(  vN-Hopf\right)  },
\end{equation}
because $\widehat{R}_{PP}^{\left(  vN-Hopf\right)  }$, $\widehat{R}%
_{PP}^{\left(  vN-Hopf\right)  \dagger}$ and
$\widehat{R}_{QQ}^{\left( vN-Hopf\right)  }$,
$\widehat{R}_{QQ}^{\left(  vN-Hopf\right)  \dagger}$ are mutually
orthogonal.
\end{proof}

\section{Conclusion}

Thus, we have introduced a couple of new bialgebras derived from
$U_{q}\left( sl_{2}\right)  $ which contain idempotents (hence some
zero divisors). In some special cases explicit formulas for
$R$-matrices are presented. We define near-$R$-matrices which
satisfy the von Neumann regularity condition.

In a similar way one can consider an analog of $U_{q}\left(
sl_{n}\right)  $ furnished by a suitable and more cumbersome family
of idempotents. Also, it would be worthwhile to investigate
supersymmetric versions of the presented structures.

Hopefully, this approach will be able to facilitate a further
research of bialgebras splitting into direct sums, which is a new
way of generalizing the standard Drinfeld-Jimbo algebras.

\begin{acknowledgements}
One of the authors (S.D.) is thankful to J.~Cuntz, P.~Etingof,
L.~Kauffman, U.~Kr\"{a}hmer, G.~Ch.~Kurinnoj, B.~V.~Novikov,
J.~Okninski, S.~A.~Ovsien\-ko, D.~Radford, C.~Ringel, J.~Stasheff,
E.~Taft, T.~Timmermann, S.~L.~Wo\-ro\-nowicz for numerous and
helpful discussions, also he is grateful to the Alexander von
Humboldt Foundation for valuable support and to M. Zirnbauer for
kind hospitality at the Institute of Theoretical Physics, Cologne
University, where this paper was finished. Both authors are indebted
to L.\ L. Vaksman\footnote{Memorial Page:
\texttt{http://webusers.physics.umn.edu/\~{}duplij/vaksman}} for
stimulating communications related to the structure of quantum
universal enveloping algebras. \end{acknowledgements}

\end{document}